\documentclass[12pt,leqno]{article}
\usepackage{amsmath, amsthm, amsfonts, amssymb, color}
\usepackage{mathrsfs}
\theoremstyle{definition}

\newcommand{\scr}[1]{\mathscr #1}
\definecolor{wco}{rgb}{0.5,0.2,0.3}

\numberwithin{equation}{section} \theoremstyle{remark}

\newcommand{\ua}{\uparrow}

\title{
{\bf Hypercontractivity and Applications for Stochastic Hamiltonian Systems}
\footnote{Supported in part by NNSFC(11431014) and Start-Up Fund of Tianjin University.}
}
\author{
{\bf Feng-Yu Wang  }\\
  \footnotesize{Center of Applied Mathematics, Tianjin University, Tianjin 300072, China}\\
    \footnotesize{Department of Mathematics, Swansea University, Singleton Park, SA2 8PP, UK}\\
\footnotesize{Email: \tttext{wangfy@bnu.edu.cn}; \tttext{F.Y.Wang@swansea.ac.uk}} }

\begin{document}
\def\tttext#1{{\normalfont\ttfamily#1}}
\def\R{\mathbb R}  \def\ff{\frac} \def\ss{\sqrt} \def\B{\mathbf B}
\def\N{\mathbb N} \def\kk{\kappa} \def\m{{\bf m}}
\def\dd{\delta} \def\DD{\Delta} \def\vv{\varepsilon} \def\rr{\rho}
\def\<{\langle} \def\>{\rangle} \def\GG{\Gamma} \def\gg{\gamma}
  \def\nn{\nabla} \def\pp{\partial} \def\EE{\scr E}
\def\d{\text{\rm{d}}} \def\bb{\beta} \def\aa{\alpha} \def\D{\scr D}
  \def\si{\sigma} \def\ess{\text{\rm{ess}}}
\def\beg{\begin} \def\beq{\begin{equation}}  \def\F{\scr F}
\def\Ric{\text{\rm{Ric}}} \def\Hess{\text{\rm{Hess}}}
\def\e{\text{\rm{e}}} \def\ua{\underline a} \def\OO{\Omega}  \def\oo{\omega}
 \def\tt{\tilde} \def\Ric{\text{\rm{Ric}}}
\def\cut{\text{\rm{cut}}} \def\P{\mathbb P}
\def\C{\scr C}     \def\E{\mathbb E}
\def\Z{\mathbb Z} \def\II{\mathbb I}
  \def\Q{\mathbb Q}  \def\LL{\Lambda}\def\L{\scr L}
  \def\B{\scr B}    \def\ll{\lambda}
\def\vp{\varphi}\def\H{\mathbb H}\def\ee{\mathbf e}

\maketitle
\begin{abstract} The hypercontractivity
is proved  for the Markov semigroup associated with a class of   stochastic Hamiltonian systems on Hilbert spaces.
 Consequently,  the  Markov semigroup converges exponentially  to the invariant probability measure  in  entropy and is compact for large  time.
These strengthen the hypocoercivity results derived in the literature. Since the log-Sobolev inequality is invalid, we introduce a new argument to prove the hypercontractivity using coupling and
 dimension-free Harnack inequality.   The main results are illustrated by  concrete examples of the kinetic Fokker-Planck equation and highly degenerate diffusion processes.
\end{abstract} \noindent

 \vskip 2cm

\section{Introduction}

To motivate the present study, we first recall the famous hypocoercivity result of C. Villani \cite{V}. Consider the following degenerate SDE (stochastic differential equation) for $(X_t,Y_t)$ on $\R^d\times\R^d$:
\beq\label{V1}\beg{cases} \d X_t= Y_t\,\d t,\\
\d Y_t= \{\nn V(X_t)-Y_t\}\d t+\ss 2\,\d W_t,\end{cases}\end{equation}
where $V\in C^2(\R^d)$   such that
$$\mu(\d x,\d y):=\e^{V(x)-\ff 1 2 |y|^2}\d x\d y$$ is a probability measure on $\R^d\times\R^d$,  and $W_t$ is the $d$-dimensional Brownian motion. This type   degenerate SDE is known as $``$Stochastic Hamiltonian System (Abbrev. SHS)''  in probability theory (see \cite{Wu}),  and the distribution density of the solution solves the kinetic Fokker-Planck equation (see \cite{V}). Let $P_t$ be the Markov semigroup for the solution of \eqref{V1}. According to \cite[Theorem 35]{V}, if there exists a constant $C>0$ such that
$$ |\nn^2 V|\le C(1+|\nn V|)$$ and the following Poincar\'e inequality holds for $\mu_1(\d x):= \mu(\d x\times\R^d)$:
$$ \mu_1(f^2)\le C\mu_1(|\nn f|^2),\ \ f\in C_b^1(\R^d), \mu_1(f)=0,$$
then for some constants $c,\ll>0$ one has
\beq\label{PV3} \mu(|\nn P_tf|^2+(P_tf)^2)\le c\e^{-\ll t}\mu(|\nn f|^2+f^2),\ \ f\in C_b^1(\R^{2d}), \mu(f)=0, t\ge 0.\end{equation}
See \cite{DMS,D,GM,GS} and references within for $L^2$-exponential convergence of  the same type degenerate diffusion semigroups.     The methodology used in these papers relies heavily on the explicit formulation of the invariant probability measure $\mu$. In this paper, we investigate  the hypercontractivity, a stronger property than the $L^2$-exponential convergence, for more general degenerate diffusion processes with  inexplicit  invariant probability measures. 

The model we investigate here is the following SHS on $\H:=\H_1\times\H_2$, where $\H_1$ and $\H_2$ are two  separable Hilbert spaces:
\beq\label{1.1}\beg{cases} \d X_t= (AX_t+BY_t)\,\d t,\\
\d Y_t= Z(X_t,Y_t)\d t +\si\d W_t,\end{cases}\end{equation}
where \beg{enumerate} \item[$\bullet$] $A$ is a densely defined (possibly unbounded) linear operator on $\H_1$;
 \item[$\bullet$] $B$ is a bounded linear operator from $\H_2$ to $\H_1$;
\item[$\bullet$] $Z$ is a densely defined map from $\H$ to $\H_2$;
\item[$\bullet$]   $\si$ is a linear operator on $\H_2$;    \item[$\bullet$] $W_t$ is the cylindrical Brownian motion on $\H_2$, i.e.
$$W_t= \sum_{i\ge 1} B_t^i e_i$$ for independent one-dimensional Brownian motions $\{B_t^i\}_{i\ge 1}$ and orthonormal basis $\{e_i\}_{i\ge 1}$ of $\H_2$. \end{enumerate}

See \cite{GW, WZ, WZ2} for results on the existence and uniqueness of (mild) solutions, as well as  Harnack inequality  and gradient estimate of the associated Markov semigroup $P_t$. We intend to find out   explicit conditions ensuring the existence and uniqueness of the invariant probability measure $\mu$ (whose formulation is in general unknown)  and, furthermore,  the hypercontractivity of $P_t$.

According to Nelson \cite{N}, $P_t$ is called hypercontractive if it has an invariant probability measure $\mu$ such that 
$$\|P_t\|_{L^2(\mu)\to L^4(\mu)}:=\sup\{\|P_tf\|_{L^4(\mu)}:\ \mu(f^2)\le 1\}=1\ \text{for\ some\ } t>0.$$  By the semigroup property and the interpolation theorem, the norm $\|\cdot\|_{L^2(\mu)\to L^4(\mu)}$ can be replaced by
$\|\cdot\|_{L^p(\mu)\to L^q(\mu)}$ for any $(p,q)\in (1,\infty)$ with $q>p$. As applications of the hypercontractivity, we will prove the compactness of $P_t$  for large $t>0$ and the exponential convergence in entropy.  %$$\mu((P_tf)\log P_t f)\le c\e^{-\ll t}\mu(f\log f),\ \ t\ge 0, f\ge 0, \mu(f\log f)=1$$ for some constants $c,\ll>0$.

Due to   L. Gross (see e.g. \cite{G}), the hypercontractivity of $P_t$ follows from    the log-Sobolev inequality
 $$\mu(f^2\log f^2)-\mu(f^2)\log \mu(f^2)\le C \EE(f,f),\ \ f\in \D(\EE)$$ for some constant $C>0$, where $(\EE,\D(\EE))$ is the associated energy form. Because of this result, the log-Sobolev inequality   has been intensively investigated for  forty years.   However,  since the energy form $\EE$ associated with  \eqref{1.1} satisfies
  $$\EE(f,f)=  \mu(|\si^* \nn_y f|^2) =0$$ for $f\in C_b^1(\H)$ with $f(x,y)$ depending only on $x$,  
   the log-Sobolev inequality is invalid. So, to prove the hypercontractivity we   need to develop a new argument.

 The remainder of the paper is organized as follows. In Section 2, we introduce a general result on the hypercontractivity using coupling and dimension-free Harnack inequality initiated from \cite{W97}. This result is then applied in Sections 3 and 4 to finite- and infinite-dimensional SHS respectively.   Finally, concrete examples  are presented in Section 5 to illustrate our main results.

 \section{Hypercontractivity using Harnack inequality}

 In this section,   we introduce a general result on the hypercontractivity using Harnack inequality.
 The basic idea of the study goes back to \cite{W97} for elliptic diffusion  semigroups on manifolds, see also \cite{BWY} for a recent study of functional SDEs.

 For  a probability space $(E,\B,\mu)$,  let $P_t$ be a Markov semigroup on $\B_b(E)$ such that $\mu$ is $P_t$-invariant, i.e. $\mu(P_tf)=\mu(f)$ for $f\in L^1(\mu)$ and $t\ge 0$.
 Recall that   a process $(X_t, Y_t)$ on $E\times E$ is called a coupling of the Markov process with   semigroup $P_t$, if
  $$(P_t f)(X_0)= \E \big(f(X_t)|X_0\big),\ \ (P_t f)(Y_0)= E\big(f(Y_t)|Y_0\big),\ \ f\in \B_b(E), t\ge 0.$$

 \beg{thm}\label{T2.1} Assume that the following three conditions  hold for some    measurable functions $\rr: E\times E\to (0,\infty)$  and
 $\phi: [0,\infty)\to (0,\infty)$ with $\lim_{t\to\infty} \phi(t)=0$:
 \beg{enumerate} \item[$(i)$] There exists two constants $t_0,c_0>0$ such that
 $$(P_{t_0}f(\xi))^2\le (P_{t_0}f^2(\eta)) \e^{c_0 \rr(\xi,\eta)^2},\ \ f\in \B_b(E), \xi,\eta\in E;$$
 \item[$(ii)$] For any $(X_0,Y_0)\in E\times E$, there exists a coupling $(X_t,Y_t)$
 associated to $P_t$ such that $$\rr(X_t,Y_t)\le \phi (t)\rr(X_0,Y_0),\ \ t\ge 0;$$
 \item[$(iii)$] There exists $\vv>0$ such that $(\mu\times\mu)(\e^{\vv \rr^2})<\infty.$ \end{enumerate}
 Then $\mu$ is the unique invariant probability measure and $P_t$ is hypercontractive.  
  Consequently, $P_t$ is  compact in $L^2(\mu)$ for large $t>0$, and there exist  constants $c,\ll>0$ such that
 \beq\label{ENT} \beg{split} &\mu((P_tf)\log P_t f)\le c\e^{-\ll t} \mu(f\log f),\ \ t\ge 0, f\ge 0, \mu(f)=1;\\
 &\|P_t f -\mu(f)\|_{L^2(\mu)}\le c\e^{-\ll t} \|f-\mu(f)\|_{L^2(\mu)},\ \ f\in L^2(\mu), t\ge 0.\end{split}\end{equation}  \end{thm}

     To prove this result, we introduce two propositions on the hypercontractivity and applications for bounded linear operators. The first is generalized from 
   \cite{W04} where symmetric Markov operators are considered.

\beg{prp}\label{P} Let $P$ be a bounded  linear operator on $L^2(\mu)$ such that $P1=1$ and $\mu$ is $P$-invariant, i.e. $\mu(P f)=\mu(f)$ for $f\in L^2(\mu).$
If $\|P\|_{L^2(\mu)\to L^4(\mu)}^4<2$, then
\beg{enumerate} \item[$(1)$] $\|P-\mu\|_{L^2(\mu)}:= \sup \{\|Pf-\mu(f)\|_{L^2(\mu)}:\ \mu(f^2)\le 1\}<1;$
\item[$(2)$]  $\|P^n\|_{L^2(\mu)\to L^4(\mu)} =1$    for large enough $n\in \mathbb N.$ \end{enumerate} \end{prp}

\beg{proof} (1)  Let $\dd(P):= \|P\|_{L^2(\mu)\to L^4(\mu)}^4<2.$   For any
 $f\in L^2(\mu)$ with $\mu(f^2)=1$
and $\mu(f)=0,$ we intend to prove
\beq\label{*8} \mu((Pf)^2)\le \inf_{\vv\in (0,1)} \ff{\ss{8\vv^2 + \dd(P)} -3\vv}{1-\vv}.\end{equation}
Without loss of generality, we assume
$\mu((Pf)^3) \ge 0,$ otherwise it suffices to replace $f$ by   $-f$.
 For any $\vv\in (0,1),$  let $g_\vv =\ss{\vv} +\ss{1-\vv}f$. Then $\mu(g_\vv^2)=1.$
Since $P1 =1,$ $\mu(Pf)   =\mu(f)=0, \mu((Pf)^3)\ge 0, \mu(g_\vv^2)=1$ and   $\mu((Pf)^4)\ge
\mu((Pf)^2)^2,$ we have
\beg{equation*}\beg{split} &\dd(P)
 \ge \mu((Pg_\vv)^4)\\
 & = \vv^2 +(1-\vv)^2 \mu((Pf)^4) +6\vv(1-\vv) \mu((Pf)^2)
 + 4\vv^{\ff 3 2 }\ss{1-\vv}\mu(Pf)
+4\ss\vv (1-\vv)^{\ff 3 2}\mu((Pf)^3)\\
&\ge (1-\vv)^2 \mu((Pf)^2)^2 +6\vv(1-\vv) \mu((Pf)^2)
+\vv^2.\end{split}
\end{equation*} This implies \eqref{*8}.  According to the calculations in \cite[pages 2632-2633]{W04},   $\dd(P)<2$ and \eqref{*8} imply
$$\|P-\mu\|_{L^2(\mu)}^2 \le  \inf_{\vv\in
(0,1)}\ff{\ss{8\vv^2 + \dd(P)} -3\vv}{1-\vv}<1. $$

(2)  For $f\in L^2(\mu)$ with $\mu(f^2)=1,$ let $\hat f= f-\mu(f).$ We have $\mu(P^m\hat f)=0, m\ge 1$. Let  $\theta:= \|P-\mu\|_{L^2(\mu)}$. Then
 $$ \mu((P^m\hat f)^2)\le \theta^{2m}\mu({\hat f}^2),\ \ m\ge 1,$$ so that
\beg{equation*}\beg{split} &\mu((P^{m+1}f)^4)  = \mu(f)^4 + 4 \mu(f) \mu((P^{m+1}\hat f)^3) + 6 \mu(f)^2 \mu((P^{m+1}\hat f)^2) +\mu((P^{m+1}\hat f)^4)\\
&\le \mu(f)^4 + 4\|P\|_{L^2(\mu)\to L^3(\mu)}^3|\mu(f)| \mu((P^{m}\hat f)^2)^{\ff 3 2} \\
&\quad + 6 \mu(f)^2 \mu((P^{m+1}\hat f)^2)
+\|P\|_{L^2(\mu)\to L^4(\mu)}^4\mu((P^{m}\hat f)^2)^2\\
&\le \mu(f)^4 + 4 \|P\|_{L^2(\mu)\to L^3(\mu)}^3 \theta^{3m}|\mu(f)| \mu({\hat f}^2)^{\ff 3 2}\\
&\quad + 6 \theta^{2(m+1)}\mu(f)^2 \mu({\hat f}^2) +\|P\|_{L^2(\mu)\to L^4(\mu)}^4\theta^{4m}\mu({\hat f}^2)^2.\end{split}\end{equation*} Since   $\theta\in (0,1) $ due to (1),    $\|P\|_{L^2(\mu)\to L^3(\mu)}\le \|P\|_{L^2(\mu)\to L^4(\mu)}<\infty$, and
$$2|\mu(f)|\mu({\hat f}^2)^{\ff 3 2} \le \mu(f)^2 \mu({\hat f}^2)+\mu({\hat f}^2)^2,$$ this implies that for large enough $m\ge 1,$
$$\mu((P^{m+1}f)^4) \le \mu(f)^4 + 2 \mu(f)^2 \mu({\hat f}^2) +\mu({\hat f}^2)^2 =\mu(f^2)^2=1.$$ Therefore,  $\|P^n\|_{L^2(\mu)\to L^4(\mu)}\le 1$ holds for large enough $n\ge 1.$
\end{proof}

Next, we present a result on exponential convergence  implied by the hypercontractivity, which is well known in the literature of symmetric Markov semigroups.

 \beg{prp}\label{PP2.3} Let $P$ be a posivity-preserving   linear operator on $L^1(\mu)$ such that $\mu$ is $P$-invariant  and $\|P\|_{L^p(\mu)\to L^q(\mu)}\le 1$ holds for some constants $q>p>1.$ Then
 \beq\label{ET}\mu((Pf)\log Pf)\le \ff{(p-1)q}{p(q-1)} \mu(f\log f),\ \ \ f\ge 0, \mu(f)=1.\end{equation} Consequently,
 \beq\label{ET'} \mu\big((Pf)^2\big) \le \ff{(p-1)q}{p(q-1)} \mu(f^2),\ \ \ f\in L^2(\mu), \mu(f)=0.\end{equation} \end{prp}

 \beg{proof} Let $f\in L^2(\mu)$ with $\mu(f)=0$.  By applying \eqref{ET} to $f_s:=\ff{1+sf}{1+s\mu(f)}$, multiplying with $s^{-2}$ and letting $s\to 0$, we prove \eqref{ET'}. So, it suffices to prove \eqref{ET}.
 For any $\vv\in (0,p-1)$, let
 $$r= \ff{p-1-\vv}{(1+\vv)(p-1)},\ \ \ \dd(\vv)= \ff{p(q-1)\vv}{(p-1-\vv) q+\vv p}.$$ Then
 $$\ff 1 {1+\vv} = r +\ff{1-r} p,\ \ \ff 1 {1+\dd(\vv)}= r+ \ff{1-r} q.$$ Since $\|P\|_{L^1(\mu)}=1$ and $\|P\|_{L^p(\mu)\to L^q(\mu)}\le 1$,   Riesz-Thorin's interpolation theorem implies
  $\|P\|_{L^{1+\vv}(\mu)\to L^{1+\dd(\vv)}(\mu)}\le 1.$ So, for any $f\in \B_b^+(E)$ with $\mu(f)=1,$
 $$\int_E (Pf^{\ff 1{1+\vv}})^{1+\dd(\vv)} \d\mu \le 1,\ \ \vv\in (0,p-1).$$ Since the equality holds for $\vv=0$, this implies
 $$\ff{\d}{\d\vv}\Big|_{\vv=0} \int_E (Pf^{\ff 1{1+\vv}})^{1+\dd(\vv)} \d\mu\le 0,$$  which is equivalent to \eqref{ET}.
 \end{proof}

 \beg{proof}[Proof of Theorem \ref{T2.1}] (a)  According to \cite[Proposition 3.1]{WY11}, $(i)$ implies that $\mu$ is the unique invariant probability measure of $P_{t_0}$, and $P_{t_0}$  has a  density with respect to $\mu$. So, by \cite[Theorem 2.3]{Wu}, if $\|P_t\|_{L^2(\mu)\to L^4(\mu)}<\infty$ then $P_{t_0+t}$ is compact in $L^2(\mu)$.
  Therefore, according to Propositions \ref{P} and \ref{PP2.3}, it remains to prove   $\|P_{t}\|_{L^2(\mu)\to L^4(\mu)}^4<2$ for large enough $t>0.$

(b) Let $f\in \B_b(E)$ with $\mu(f^2)\le 1.$ By $(i)$ and $(ii)$  we have
\beg{align*} &(P_{t+t_0} f(\xi))^2 \le \E (P_{t_0} f(X_t))^2 \le \E\Big[(P_{t_0} f^2(Y_t))\e^{c_0\rr(X_t,Y_t)^2}\Big]\\
&\le (P_{t_0+t} f^2(\eta)) \e^{c_0\phi(t)^2 \rr(\xi,\eta)^2},\ \ t\ge 0, (\xi,\eta)\in E\times E. \end{align*}  Equivalently,
$$(P_{t_0+t}f(\xi))^2 \e^{-c_0\phi(t)^2\rr(\xi,\eta)^2} \le P_{t_0+t} f^2(\eta),\ \ t\ge 0, (\xi,\eta)\in E\times E.$$ Integrating with respect to $\mu(\d \eta)$ gives
$$(P_{t_0+t} f(\xi))^2 \int_E \e^{-c_0 \phi(t)^2 \rr(\xi,\eta)^2}\mu(\d\eta) \le \int_E P_{t_0+t}f^2(\eta)\mu(\d \eta) =\mu(f^2)\le 1,\ \ t\ge 0,\xi\in E.$$ Thus,
$$(P_{t_0+t} f(\xi))^4\le \ff 1 { \big(\int_E \exp[-c_0 \phi(t)^2 \rr(\xi,\eta)^2]\mu(\d\eta)\big)^2},\ \ \ \mu(f^2)\le 1, t\ge 0, \xi\in E.$$ Then
by Jensen's inequality, for $t\ge 0$
\beq\label{BG}\beg{split} &\sup_{\mu(f^2)\le 1} \int_E (P_{t+t_0}f(\xi))^4\mu(\d \xi)  \le \int_E \ff{\mu(\d\xi)}{(\int_E \exp[-c_0\phi(t)^2\rr(\xi,\eta)^2]\mu(\d\eta))^2}\\
&\le \int_E  \bigg(\int_E\e^{c_0\phi(t)^2\rr(\xi,\eta)^2}\mu(\d\eta)\bigg)^2\mu(\d\xi)
 \le \int_{E\times E} \e^{2c_0\phi(t)^2\rr(\xi,\eta)^2}\mu(\d\xi)\mu(\d\eta).\end{split}\end{equation}  Since   $\lim_{t\to \infty} \phi(t)=0$, it follows from  $(iii)$ that
  $$\lim_{t\to\infty} \int_{E\times E} \e^{2c_0\phi(t)^2\rr(\xi,\eta)^2}\mu(\d\xi)\mu(\d\eta)=1.$$  Combining this with \eqref{BG} we prove $\|P_t\|_{2\to 4}^4<2$ for large enough $t>0.$ \end{proof}

 \section{ Hypercontractivity for finite-dimensional SHS}

 In this section, we consider the equation \eqref{1.1} with $\H=\R^{m+d}$ for some $m,d\ge 1.$   Let $\|\cdot\|$ denote the operator norm.
  To verify conditions $(i)$-$(iii)$ in Theorem \ref{T2.1}, we make the following assumptions.

 \beg{enumerate} \item[$(A1)$] $\si$ is invertible and Rank$[B, AB,\cdots, A^{m-1} B]=m.$
 \item[$(A2)$] $Z: \R^{m+d}\to \R^d$ is Lipschitz continuous.
 \item[$(A3)$] There exist constants $r,\theta>0$ and  $r_0\in (-\|B\|^{-1}, \|B\|^{-1})$    such that
 \beg{equation*}\beg{split} &\big\<r^2(x-\bar x)+ rr_0B(y-\bar y), A(x-\bar x) +B(y-\bar y)\big\>\\
 & +\big\<Z(x,y)-Z(\bar x,\bar y), y-\bar y +rr_0 B^* (x-\bar x)\big\>\\
  &\le -\theta (|x-\bar x|^2 +|y-\bar y|^2),\ \ \ \ (x,y), (\bar x,\bar y)\in \R^{m+d}.\end{split}\end{equation*}
 \end{enumerate}

The rank condition in $(A1)$ is known as Kalman's  condition, when $\si$ is invertible it   is equivalent to the H\"ormander condition. We will prove the Harnack inequality in condition $(i)$ using $(A1)$ and $(A2)$, and verify  conditions $(ii)$ and $(iii)$ by   Assumption $(A3)$.

 \beg{thm}\label{T1.1} Assume $(A1)$, $(A2)$ and $(A3)$. Let $P_t$ be the Markov semigroup associated with $\eqref{1.1}$. Then
 \beg{enumerate} \item[$(1)$] $P_t$ has a unique invariant probability measure $\mu$ and $\mu(\e^{\vv |\cdot|^2})<\infty$ for some $\vv>0;$
 \item[$(2)$] $P_t$ is hypercontractive, i.e. $\|P_t\|_{2\to 4}=1$ for large $t>0;$
  \item[$(3)$]   $P_t$ is compact in $L^2(\mu)$ for large $t>0,$ and there exist constants $c,\ll>0$ such that $\eqref{ENT}$ holds.
 \end{enumerate}\end{thm}

In a similar spirit of \eqref{PV3},   under  a generalized curvature condition \cite{B} proved the following entropy-information inequality for some constants $c,\ll>0$:
$$\mu\big((P_tf)\log P_t f+ (P_tf)|\nn \log P_t f|^2\big)\le c\e^{-\ll t} \mu\big(f\log f+f|\nn \log f|^2\big),\ \ f\ge 0, \mu(f)=1, t\ge 0.$$
This does not imply the entropy inequality in \eqref{ENT}.

\

 According to Theorem \ref{T2.1} and Proposition \ref{PP2.3},  Theorem \ref{T1.1} follows from  the following three lemmas which correspond   to conditions $(i)$-$(iii)$ respectively.
The first lemma provides the desired Harnack inequality. Although the Harnack  inequality has been  investigated in \cite{GW,WZ} for SHS,
the resulting results are not enough for our   purpose:  the   inequality established in \cite{GW} (see Corollary 4.2 therein) contains a worse exponential term, while the assumption (H) in \cite{WZ} does not hold if $Z$ is not second order differentiable. So, we present below a new version  of   Harnack inequality for SHS using coupling by change of measures. See
\cite[Chapter 1]{Wbook} for more results on the coupling by change measures and applications.

 \beg{lem}\label{L1} Assume $(A1)$ and $(A2)$. For any $t_0>0$,  there exists a constant $c_0>0$ such that
 $$(P_{t_0}f)^2(\xi)\le (P_{t_0}f^2(\eta))\e^{c_0 |\xi-\eta|^2},\ \ f\in \B_b(\R^{m+d}), \xi,\eta \in \R^{m+d}.$$\end{lem}

 \beg{proof} Let $(X_t,Y_t)$ solve the equation \eqref{1.1} with $(X_0, Y_0)=\eta\in\R^{m+d}$, and let $(\bar X_t, \bar Y_t)$ solve the following equation with $(\bar X_0, \bar Y_0)=\xi\in\R^{m+d}$:
 \beq\label{E2} \beg{cases} \d \bar X_t= (A\bar X_t+ B\bar Y_t)\,\d t,\\
 \d \bar Y_t= \Big\{Z(X_t,Y_t)  + \dfrac {Y_0-\bar Y_0} {t_0}   + \dfrac \d{\d t}\big(t(t_0-t)B^*\e^{(t_0-t)A^*}b\big)\Big\}\d t +\si \d W_t,\end{cases}\end{equation}
 where $b\in\R^m$ is to be determined   such that $(X_{t_0},Y_{t_0}) =(\bar X_{t_0},\bar Y_{t_0}).$ It is easy to see that
 $$\beg{cases} \ff{\d}{\d t} (X_t-\bar X_t)=  A(X_t-\bar X_t)+ B(Y_t-\bar Y_t),\\
 \ff{\d}{\d t} (Y_t- \bar Y_t)=   \ff 1 {t_0} (\bar Y_0-  Y_0)  - \ff{\d}{\d t}\big\{t(t_0-t)B^*\e^{(t_0-t)A^*}b\big\}.\end{cases}$$ Then
 \beq\label{E3} Y_t-\bar Y_t= \ff{t_0-t}{t_0}(Y_0-\bar Y_0)  - t(t_0-t) B^* \e^{(t_0-t)A^*}b,\end{equation} and
\beq\label{E4}\beg{split} X_t-\bar X_t &= \e^{At} (X_0-\bar X_0)+\int_0^t\e^{A(t-s)}B(Y_s-\bar Y_s)\d s\\
&= \e^{At} (X_0-\bar X_0)+ \bigg(\int_0^t \e^{A(t-s)}\ff{t_0-s}{t_0} \d s\bigg)B(Y_0-\bar Y_0)\\
&\quad - \bigg(\int_0^t s(t_0-s)\e^{A(t-s)}BB^*\e^{(t_0-s)A^*}\d s\bigg)b.\end{split}\end{equation} We now    take
\beq\label{E5} b=Q_{t_0}^{-1} \bigg\{\e^{t_0 A}(X_0-\bar X_0)+ \bigg(\int_0^{t_0}\ff{t_0-s}{t_0} \e^{A(t_0-s)}\d s\bigg)B(Y_0-\bar Y_0)\bigg\},\end{equation}
where, according to \cite[\S 3]{S}, the rank condition in $(A1)$ ensures the invertibility of the $m\times m$-matrix
$$Q_{t_0}:= \int_0^{t_0} s(t_0-s)\e^{A(t_0-s)}BB^*\e^{(t_0-s)A^*}\d s,$$  see (1) in the proof of \cite[Theorem 4.2]{WZ} for details.
Then \eqref{E3}-\eqref{E5} imply $(X_{t_0}, Y_{t_0})= (\bar X_{t_0}, \bar Y_{t_0}).$

 In order to establish the Harnack inequality using Girsanov's theorem, let
$$\psi_t=Z(X_t,Y_t)-Z(\bar X_t,\bar Y_t) +\ff 1 {t_0} (Y_0-\bar Y_0)+\ff{\d}{\d t} \big\{t(t_0-t)B^*\e^{(t_0-t)A^*}b\big\},\ \ t\in [0,t_0].$$ Since $Z$ is Lipschitz continuous,  \eqref{E3}, \eqref{E4} and \eqref{E5} imply
\beq\label{E6} |\psi_t|^2 \le c_1 (|X_0-\bar X_0|^2+|Y_0-\bar Y_0|^2)= c_1|\xi-\eta|^2,\ \ t\in [0,t_0]\end{equation} for some constant $c_1>0.$ Moreover, according to the definition of $\psi$, \eqref{E2} can be reformulated as
$$\beg{cases} \d\bar X_t=  (A \bar X_t + B \bar Y_t)\,\d t,\\
 \d\bar Y_t =  Z(\bar X_t,\bar Y_t)\d t +\si \d \bar W_t,\end{cases}$$ where
 $$\bar W_t:= W_t +\si^{-1} \int_0^t\psi_s\d s,\ \ t\in [0,t_0].$$ Let
 \beq\label{E7} R:= \exp\bigg[-\int_0^{t_0} \<\si^{-1}\psi_t,\d W_t\>-\ff 1 2 \int_0^{t_0} |\si^{-1}\psi_t|^2 \d t\bigg].\end{equation} By \eqref{E6} and Girsanov's theorem, $\tt W_t$ is a $d$-dimensional Brownian motion under the probability measure $\d\Q:= R\d\P$.  Therefore, by the weak uniqueness of the equation \eqref{1.1}  and using $(X_{t_0}, Y_{t_0})= (\bar X_{t_0}, \bar Y_{t_0})$, we obtain
\beg{equation*}\beg{split} (P_{t_0}f(\xi) )^2 &= \big(\E[Rf(\bar X_{t_0},\bar Y_{t_0})]\big)^2 = \big(\E[Rf(X_{t_0},Y_{t_0})]\big)^2\\
&\le (\E R^2) \E f^2 (X_{t_0}, Y_{t_0}) = (P_{t_0} f^2 (\eta)) \E R^2.\end{split}\end{equation*} Noting that   \eqref{E6} and \eqref{E7} imply $\E R^2\le \e^{c_0 |\xi-\eta|^2}$ for some constant $c_0>0,$ we finish the proof.
 \end{proof}

\beg{lem}\label{L2} If $(A3)$ holds, then there exist two constants $c,\ll>0$ such that for any two solutions $(X_t,Y_t)$ and $(\tt X_t,\tt Y_t)$  of $\eqref{1.1}$,
$$ |X_t-\tt X_t|^2 +|Y_t-\tt Y_t|^2 \le c \e^{-\ll t} (|X_0-\tt X_0|^2 +|Y_0-\tt Y_0|^2),\ \ t\ge 0.$$\end{lem}

\beg{proof} Obviously, $X_t-\tt X_t$ solves the ODE
\beq\label{D0}\beg{cases} \ff{\d}{\d t} (X_t-\tt X_t)=  A(X_t-\tt X_t)+ B(Y_t-\tt Y_t),\\
 \ff{\d}{\d t} (Y_t- \tt Y_t)=   \big(Z(X_t,Y_t)-Z(\tt X_t,\tt Y_t)\big)\d t.\end{cases}\end{equation} Since $r_0\in (-\|B\|^{-1}, \|B\|^{-1})$, for any $r>0$  there exists a constant $C>1$ such that
 \beq\label{D1}\beg{split}  &\ff 1 C (|X_t-\tt X_t|^2+|Y_t-\tt Y_t|^2) \\
 &\le \Phi_t  := \ff {r^2} 2 |X_t-\tt X_t|^2 +\ff 1 2 |Y_t-\tt Y_t|^2 + r r_0 \<X_t-\tt X_t, B(Y_t-\tt Y_t)\>\\
 &\le C(|X_t-\tt X_t|^2+|Y_t-\tt Y_t|^2),\ \ t\ge 0.\end{split}\end{equation} Combining this with \eqref{D0} and $(A3)$, we obtain
 $$\d\Phi_t \le -\theta (|X_t-\tt X_t|^2 +|Y_t-\tt Y_t|^2) \le -\ff \theta C  \Phi_t\d t.$$ Therefore, $\Phi_t \le \Phi_0 \e^{-\theta t/C}.$
 This together with   \eqref{D1} implies the desired estimate. \end{proof}

 \beg{lem}\label{L3} If $(A3)$ holds, then $P_t$ has an invariant probability measure $\mu$ such that $\mu(\e^{\vv |\cdot|^2})<\infty$   for some constant $\vv>0.$ \end{lem}

 \beg{proof} Let $(X_t,Y_t)$ solve \eqref{1.1} with $(X_0,Y_0)=0\in\R^{m+d}.$ By a standard tightness argument, it suffices to prove
 \beq\label{D2} \sup_{t\ge 0} \E \e^{\vv(|X_t|^2+|Y_t|^2)}<\infty\end{equation} for some constant $\vv>0.$ Since $r_0\in (-\|B\|^{-1}, \|B\|^{-1})$, for any $r>0$ there exists a constant $C>1$ such that
 \beq\label{D3}\beg{split}  \ff 1 C (|X_t|^2+|Y_t|^2) &\le \Psi_t  := \ff {r^2} 2 |X_t|^2 +\ff 1 2 |Y_t|^2 + rr_0 \<X_t, BY_t\>\\
 &\le C(|X_t|^2+|Y_t|^2),\ \ t\ge 0.\end{split}\end{equation} Moreover,   $(A3)$ with $(\bar x,\bar y)=0$ implies
 $$\<r^2x+rr_0 By, Ax + By\> + \<Z(x,y)-Z(0,0), y+rr_0B^* x\>\le -\theta (|x|^2+|y|^2),\ \ \ (x,y)\in \R^{m+d}.$$ Then there exist constants $c_1,c_2>0$ such that
\beg{equation*}\beg{split} &\<r^2x+rr_0 By, Ax + By\> + \<Z(x,y), y+r r_0 B^* x\>\\
&\le |Z(0,0)|\cdot|y+rB^*x|-\theta (|x|^2+|y|^2)\le c_1 -c_2 (|x|^2+|y|^2),\  \ (x,y)\in \R^{m+d}.\end{split}\end{equation*}
Thus, by \eqref{1.1}, It\^o's formula and \eqref{D3}, we may find out two constants $c_3,c_4>0$ such that
\beg{equation*}\beg{split}  \d \Psi_t &\le \big(c_3-c_2 (|X_t|^2+|Y_t|^2)\big)\d t +\<Y_t+ r B^*X_t, \si \d W_t\>\\
&\le (c_3-c_4 \Psi_t)\d t +\<Y_t+ r B^*X_t, \si \d W_t\>.\end{split}\end{equation*} By It\^o's formula, for any $\vv>0$ there exists a local martingale $M_t$ such that 
 $$\d \e^{\vv \Psi_t} \le \vv\e^{\vv \Psi_t}\Big(c_3-c_4\Psi_t +\ff{\vv^2} 2 |\si^*(Y_t+rB^*X_t)|^2\Big)\d t +\d M_t.$$   Noting that \eqref{D3}  implies  $|\si^*(Y_t+rB^*X_t)|^2\le c_5 \Psi_t$ for some constant $c_5>0,$ by taking $\vv=\ff{c_4}{c_5}$ we obtain
 $$\d \e^{\vv \Psi_t} \le \vv\e^{\vv \Psi_t}\Big(c_3-\ff 1 2 c_4\Psi_t\Big)\d t +\d M_t \le (c_6-\e^{\vv\Psi_t})\d t+\d M_t$$ for some constant $c_6\ge 1.$  Since $\e^{\vv \Psi_0}=1,$ it follows that
 $$\E \e^{\vv \Psi_t}\le  c_6,\ \ t\ge 0.$$ Because of \eqref{D3}, this implies \eqref{D2} for small $\vv>0$.
   \end{proof}

\section{Hypercontractivity for infinite-dimensional SHS}

When $\H_2$ is infinite-dimensional and  $\si$ is not Hilbert-Schmidt,  $\si W_t$  is ill defined   on $\H_2$, so that
the usual  strong solution of \eqref{1.1}  does not make sense. Alternatively,  we consider the mild solution. To this end,  we reformulate \eqref{1.1}   on $\H:=\H_1\times\H_2$ as follows:
\beq\label{E1.1}\beg{cases} \d X_t= (AX_t+BY_t-L_1X_t)\,\d t,\\
\d Y_t= \{Z(X_t,Y_t)-L_2Y_t\}\d t +\si\d W_t,\end{cases}\end{equation}where  $A:\H_1\to \H_1, B: \H_2\to\H_1$ and
$\si:\H_2\to\H_2$ are bounded linear operators;    $(L_i,\D(L_i))$ is a positive definite self-adjoint operator on $\H_i, i=1,2$;
and $Z: \H \to \H_2$ is measurable.
This equation reduces to \eqref{1.1} if we regard $A-L_1$ as one operator and combine $Z(x,y)$ with $-L_2 y$. The unbounded operator $L_2$ plays a crucial role in the study of mild solutions (see \cite{DP}), while
 $L_1$ is the counterpart of $L_2$ for the first component process $X_t$, and the bounded operator $A$ stands for a perturbation of $L_1$, see $(B3)$ below.

Let
$\<\cdot,\cdot\>, |\cdot|$ and $\|\cdot\|$   denote, respectively,  the inner product, the norm and the operator norm on  a Hilbert space. Moreover,
for a linear operator $(L,\D(L))$ on a Hilbert space, and for $\ll\in \R$, we write $L\ge\ll$ if $\<f,Lf\>\ge \ll |f|^2$ holds for all $f\in \D(L).$

To prove the hypercontractivity using Theorem \ref{T2.1},  we will need the following assumptions.

\beg{enumerate} \item[$(B1)$] $\si$ is invertible,  $L_2$ has discrete spectrum with eigenbasis $\{e_i\}_{i\ge 1}$
and corresponding eigenvalues $0<\ll_1\le \ll_2\le\cdots$ including multiplicities satisfy   $\sum_{i=1}^\infty \ff 1 {\ll_i}<\infty.$
 \item[$(B2)$] There exist two constants $K_1,K_2>0$ such that
 $$|Z(x,y)-Z(\bar x, \bar y)|\le K_1|x-\bar x|+K_2|y-\bar y|,\ \ (x,y), (\bar x,\bar y)\in \H.$$
 \item[$(B3)$] $L_1-A\ge \ll_1-\dd$ for some constant $\dd\ge 0$, $BL_2=L_1B, AL_1=L_1A$, and for any $t>0$
 $$Q_t:= \int_0^t \e^{sA} BB^* \e^{sA^*}\d s $$ is an invertible operator on $\H_1$.
 \end{enumerate}
 It is well known that $(B1)$ and $(B2)$ imply the existence and uniqueness of mild solutions for \eqref{E1.1}, see \cite{DP}. Let $P_t$ be the associated Markov semigroup.

 \beg{thm}\label{T1.2} Assume $(B1)$, $(B2)$ and $(B3)$. If
 \beq\label{C11} \ll_1>\ll':= \ff 1 2 \Big(\dd+K_2 +\ss{(K_2-\dd)^2+4K_1\|B\|}\Big),\end{equation}
 then all assertions in Theorem $\ref{T1.1}$ hold. \end{thm}

As shown in the proof of Theorem \ref{T1.1}, we  need to verify conditions $(i)$-$(iii)$ in Theorem \ref{T2.1}.   Let $(X_t,Y_t)$ be a mild solution to \eqref{E1.1}. We have

\beq\label{E1.2} \beg{cases} X_t=\e^{-(L_1-A+\dd)t }X_0 + \int_0^t \e^{-(L_1-A+\dd)(t-s)} (\dd X_s+BY_s)\d s,\\
Y_t= \e^{-L_2t} Y_0+ \int_0^t \e^{-L_2(t-s)}Z(X_s,Y_s)\d s +\xi_t,\end{cases}\end{equation}
where $$\xi_t:= \int_0^t\e^{-L_2(t-s)}\si\d W_s,\ \ \ t\ge 0.$$
Due to $(B1)$, for any $T>0$, the process
$$M_t^T:= \int_0^t\e^{-L_2(T-s)}\si\d W_s,\ \ t\in [0,T]$$ is a square integrable martingale on $\H$ with quadratic variation process
$$\<M^T\>_t=\int_0^t \|\e^{-L_2(T-s)}\si\|_{HS}^2\d s \le \|\si\|^2 \sum_{i=1}^\infty \ff 1 {2\ll_i}=:\aa_0<\infty,\ \ t\in [0,T], $$
where $\|\cdot\|_{HS}$ is the Hilbert-Schmidt norm.
This implies
\beq\label{C0}  \E \exp\Big[\ff{ |M_t^T|^2}{2+\aa_0}\Big] \le C,\ \ T>0,  t\in [0,T]\end{equation} for some constant $ C>0$.   Indeed, since
$$\d |M_t^T|^2 = 2 \<M_t^T,\d M_t^T\> +\d \<M^T\>_t,\ \ t\in [0,T],$$ by  It\^o's formula, for any $r>0$  we have
\beg{equation*}\beg{split} &\d\Big\{\exp\Big[\ff{r |M_t^T|^2+1}{\<M^T\>_t+1}\Big]\Big\}= 
\exp\Big[\ff{r |M_t^T|^2+1}{\<M^T\>_t+1}\Big] \ff {2r}{\<M^T\>_t+1} \<M_t^T,\d M_t^T\>\\
 &- \exp\Big[\ff{r|M_t^T|^2+1}{\<M^T\>_t+1}\Big]\bigg\{\ff{r |M_t^T|^2 +1-r\<M^T\>_t-r-2r^2|M_t^T|^2}{(\<M^T\>_t+1)^2}\bigg\}\d\<M^T\>_t,\ \ t\in [0,T].\end{split}\end{equation*}
Since $\<M^T\>_t\le\aa_0,$   when   $r\in (0, \ff 1 {2+\aa_0}]$ the process $\exp\big[\ff{r |M_t^T|^2+1}{\<M^T\>_t+1}\big]$ for $t\in [0,T]$ is a supmartingale. In particular,   by  taking $r=\ff 1 {2+\aa_0}$  we prove \eqref{C0}.

Since $\xi_T= M_T^T$ for any $T>0$,  \eqref{C0} implies
\beq\label{CC} \sup_{t\ge 0} \E \exp\Big[\ff{|\xi_t|^2}{2+\aa_0}\Big]  \le C. \end{equation}
We are now ready to prove the following four lemmas which imply Theorem \ref{T1.2}  according to Theorem \ref{T2.1}.

\beg{lem}\label{L4.1}  Assume $(B1)$, $(B2)$ and $(B3)$. For any $t_0>0$, there exists a constant $c_0>0$ such that
 $$(P_{t_0}f)^2(\xi)\le (P_{t_0}f^2(\eta))\e^{c_0 |\xi-\eta|^2},\ \ f\in \B_b(\H), \xi,\eta \in  \H:=\H_1\times\H_2.$$  \end{lem}

\beg{proof} Let $(X_t,Y_t)$ solve \eqref{E1.1} with  $ (X_0,Y_0)=\eta$, and let $(\bar X_t,\bar Y_t)$ solve the following equation for $(\bar X_0,\bar Y_0)=\xi$:
$$ \beg{cases} \d \bar X_t=  (A \bar X_t+B\bar Y_t-L_1\bar X_t)\d t,\\
\d \bar Y_t=  \Big\{Z(X_t,Y_t)-L_2\bar Y_t +\ff 1 {t_0}\e^{-L_2t}(Y_0-\bar Y_0) +
\e^{-L_2t}\ff{\d}{\d t} \big(t(t_0-t)B^*\e^{(t_0-t)A^*}b\big)\Big\}\d t +\si\d W_t,\end{cases}$$
where $b\in\H_1$ will be determined latter such that $(X_{t_0},Y_{t_0})=(\bar X_{t_0},\bar Y_{t_0}).$  We have
$$ \beg{cases} \d  (X_t-\bar X_t)=    \big\{A (X_t-\bar X_t)+B(Y_t-\bar Y_t)-L_1(X_t-\bar X_t)\big\}\d t,\\
\d   (Y_t-\bar Y_t) = - \Big\{ L_2(Y_t-\bar Y_t) +\ff 1 {t_0}\e^{-L_2t}(Y_0-\bar Y_0) + \e^{-L_2t}\ff{\d}{\d t} \big(t(t_0-t)B^*\e^{(t_0-t)A^*}b\big)\Big\}\d t. \end{cases}$$ Then
\beq\label{Y} Y_t-\bar Y_t = \ff{t_0-t}{t_0} \e^{-L_2t}(Y_0-\bar Y_0)- t(t_0-t) \e^{-L_2t}B^*\e^{(t_0-t)A^*}b,\ \ t\in [0,t_0],\end{equation} and,  since $BL_2=L_1B, AL_1=L_1A$,
\beq\label{X}\beg{split}  X_t-\bar X_t  = & \e^{(A-L_1)t} (X_0-\bar X_0) + \int_0^t \ff{t_0-s}{t_0} \e^{(A-L_1)(t-s)}B\e^{-L_2s}(Y_0-\bar Y_0)\d s\\
&\quad - \int_0^t s(t_0-s) \e^{(A-L_1)(t-s)}B\e^{-L_2s} B^* \e^{A^*(t_0-s)}b\,\d s \\
 = &\e^{-tL_1} \bigg\{\e^{At} (X_0-\bar X_0) +   \int_0^t\ff{t_0-s}{t_0} \e^{A(t-s)}B(Y_0-\bar Y_0)\d s \\
 &\qquad \qquad-\int_0^t s(t_0-s) \e^{A(t-s)}BB^*\e^{A^*(t_0-s)}b\,\d s\bigg\}.\end{split}\end{equation}
According to $(B3)$,   the operator
$$\tt Q_{t_0}:=  \int_0^{t_0}s(t_0-s) \e^{A(t_0-s)}BB^*\e^{A^*(t_0-s)}\d s$$ is invertible on $\H_1$. So, letting
$$b= \tt Q_{t_0}^{-1}  \bigg\{\e^{At_0} (X_0-\bar X_0) + \int_0^{t_0} \ff{t_0-s}{t_0} \e^{A(t_0-s)}B(Y_0-\bar Y_0)\d s\bigg\},$$ we conclude  from \eqref{Y} and \eqref{X} that $(X_{t_0},Y_{t_0})= (\bar X_{t_0},\bar Y_{t_0}).$ Moreover,   there exists a constant $C_1>0$ such that
\beq\label{Li} |X_t-\bar X_t|+|Y_t-\bar Y_t|\le C_1(|X_0-\bar X_0|+|Y_0-\bar Y_0|),\ \ t\in [0,t_0].\end{equation} Since $A,B$ are bounded, $\si$ is reversible, and  $Z$ is Lipschitz continuous, this implies that  the process
$$\psi_t:= \si^{-1} \Big\{Z(X_t,Y_t)-Z(\bar X_t,\bar Y_t) +\ff 1 {t_0} \e^{-L_2t} (Y_0-\bar Y_0) +\e^{-L_2t}\ff{\d }{\d t} \big(t(t_0-t) B^* \e^{(t_0-t)A^*}\big)b\Big\}$$   satisfies
$$|\psi_t|^2  \le C_2 (|X_0-\bar X_0|^2+|Y_0-\bar Y_0|^2),\ \ t\in [0,t_0] $$ for some constant $C_2>0.$
By  the Girsanove theorem, 
$$\tt W_t:= W_t+ \int_0^t   \psi_s\d s,\ \ t\in [0,t_0]$$ is a cylindrical Brownian motion on $\H_2$ under the probability measure $\d\Q:=R\,\d\P$, where
$$R:=\exp\bigg[-\int_0^{t_0} \< \psi_s,\d W_s\> -\ff 1 2\int_0^{t_0} |\psi_s|^2\d s\bigg].$$ Rewrite the equation for   $(\bar X_t,\bar Y_t)$ as
$$ \beg{cases} \d \bar X_t=  (A \bar X_t+B\bar Y_t-L_1\bar X_t)\d t,\\
\d \bar Y_t=  \big\{Z(\bar X_t,\bar Y_t)-L_2\bar Y_t\big\}\d t +\si\d \tt W_t.\end{cases}$$ By  the weak uniqueness of the mild solutions to \eqref{E1.1} and $(X_{t_0},Y_{t_0})= (\bar X_{t_0},\bar Y_{t_0}),$ we obtain
$$(P_{t_0} f(\xi))^2= (\E_\Q f(\bar X_{t_0},\bar Y_{t_0}))^2 =(\E [R f(X_{t_0},Y_{t_0})])^2 \le (P_{t_0}f^2)(\eta) \E R^2\le (P_{t_0}f^2)(\eta)
\e^{c_0|\xi-\eta|^2}$$  for some constant $c_0>0$. \end{proof}

\beg{lem}\label{L4.2}  Assume $(B1)$, $(B2)$ and $(B3)$.    Let $(X_t,Y_t)$ solve $\eqref{E1.2}$ for $X_0=Y_0=0.$ If $ \ll_1>\ll',$ then there exists a constant $\vv>0$ such that
$\sup_{t\ge 0} \E \e^{\vv (|X_t|^2 +|Y_t|^2)} <\infty.$  \end{lem}

\beg{proof} By $(B2)$, there exists a constant $c>0$ such that
$$ |Z(x,y)|\le c+ K_1 |x|+ K_2 |y|,\ \ x,y\in\H.$$
Combining this with \eqref{E1.2}, and  noting that $(B1)$ and $(B3)$ imply  $L_1-A+\dd\ge \ll_1$ and $ L_2\ge \ll_1$, we obtain
\beq\label{4.2} \beg{split} &|X_t|\le  \int_0^t \e^{-\ll_1 (t-s)} (\dd  |X_s|+ \|B\|\cdot |Y_s)\d s,\\
&|Y_t|\le \int_0^t \e^{-\ll_1 (t-s)} (c+K_1|X_s|+K_2|Y_s|)\d s +|\xi_t|.\end{split}\end{equation}
By $(B2)$ and $(B3)$, we have    $$\aa:= \ff 1 {2\|B\|}\Big( \dd-K_2 +\ss{(K_2-\dd)^2 + 4 K_1 \|B\|}\Big)\in (0,\infty).$$
Obviously, the definitions  of $\aa$ and $\ll'$ in \eqref{C11} imply
\beq\label{DEE} \ll'\aa= \aa\dd+K_1,\ \ \ \aa\|B\|+K_2= \ll'.\end{equation} So,
$$(\aa\dd   + K_1) s+ (\aa\|B\|+K_2) t= \ll' (\aa s +t),\ \ s,t\ge 0.$$  Combining this with \eqref{4.2}, we obtain
\beg{equation*}\beg{split} \aa|X_t|+|Y_t|&\le \int_0^t\e^{-\ll_1(t-s)} \Big\{c+ (\aa\dd+K_1)|X_s|+ (\aa\|B\|+K_2)|Y_s|\Big\}\d s +|\xi_t|\\
&\le \ll' \int_0^t \e^{-\ll_1(t-s)}
(\aa|X_s|+ |Y_s|)\d s +|\xi_t|+\ff c{\ll_1}.\end{split}\end{equation*} By Gronwall's inequality, this implies
\beq\label{4.3} \beg{split} \aa|X_t|+|Y_t|&\le |\xi_t|+ \ff c{\ll_1} +\ll'\int_0^t \e^{-\ll(t-s)}\Big(|\xi_s|+\ff c{\ll_1}\Big)\d s\\
&\le |\xi_t| + c_1 + \ll'\int_0^t \e^{-\ll(t-s)} |\xi_s|\d s,\ \ t\ge 0 \end{split}\end{equation} for some constant $c_1>0$ and $\ll:=\ll_1-\ll'>0.$

Finally, applying  Jensen's inequality to  the probability measure $\nu(\d s):= \ll \e^{-\ll (t-s)}\d s $ on $(-\infty, t]$, we obtain
\beg{equation*}\beg{split} &\exp\bigg[\vv \bigg(\ll'\int_0^t \e^{-\ll(t-s)} |\xi_s| \d s\bigg)^2\bigg]= \exp\bigg[\ff \vv {\ll^2} \bigg(\ll'\int_{-\infty}^t 1_{[0,t]}(s) |\xi_s|\nu(\d s)\bigg)^2\bigg]\\
&\le \int_{-\infty}^t \exp\Big[\ff {\vv (\ll')^2}{\ll^2}1_{[0,t]}(s) |\xi_s|^2\Big] \nu(\d s)\\
& \le c_2 + c_2\int_0^t \e^{-\ll(t-s)} \exp\big[c_2 \vv |\xi_s|^2\big]\d s,\ \ t,\vv\ge 0\end{split}\end{equation*} for some constant $c_2>0.$
Combining this with \eqref{CC} and \eqref{4.3}, we finish the proof.\end{proof}

\beg{lem}\label{L4.3}  Assume $(B1)$, $(B2)$ and $(B3)$. If $\ll_1>\ll',$  then $P_t$ has a unique invariant probability measure $\mu$, and
$\mu(\e^{\vv |\cdot|^2})<\infty$ holds   for some constant $\vv>0.$ \end{lem}

\beg{proof} According to \cite[Proposition 3.1]{WY11},  the Harnack inequality in Lemma \ref{L4.1} implies that $P_t$ has at most one invariant probability measure. So, it suffices to prove the existence of $\mu$ with $\mu(\e^{\vv |\cdot|^2})<\infty$    for some constant $\vv>0.$

Let $(X_t,Y_t)_{t\ge 0}$ solve \eqref{E1.1} for $X_0=Y_0=0.$ For every $t\ge 0$, let  $\mu_t$ be the distribution of $(X_t,Y_t)$, which is a probability measure on $\H$. By the Markov property, if $\mu_t$ converges weakly to a probability measure $\mu$ as $t\to\infty$, then $\mu$ is an invariant probability measure of $P_t$  and,  by Lemma \ref{L4.2} and Fatou's lemma, $\mu(\e^{\vv |\cdot|^2})<\infty$ holds   for some constant $\vv>0.$ Therefore, it remains to prove the weak convergence of $\mu_t$ as $t\to\infty$.

Consider the $L^1$-Wasserstein distance
$$W(\nu_1,\nu_2):= \inf_{\pi\in \scr C(\nu_1,\nu_2)} \int_{\H\times\H} |\cdot|\d\pi$$ for two probability measures $\nu_1$ and $\nu_2$ on $\H\times\H$, where $ \scr C(\nu_1,\nu_2)$ is the set of all couplings of these two measures. If $\mu_t$ is a $W$-Cauchy family as $t\to\infty$,
i.e.
\beq\label{CA}\lim_{t_1,t_2\to\infty} W(\mu_{t_1},\mu_{t_2})=0,\end{equation} then it converges weakly as $t\to\infty$, see e.g. \cite[Theorem 5.4 and Theorem 5.6]{Chen}.

To prove \eqref{CA}, for any  $t_2>t_1>0$, let $(X_t,Y_t)_{t\ge 0}$ solve \eqref{E1.1} for $X_0=Y_0=0$, and let $(\tt X_t,\tt Y_t)_{t\ge t_2-t_1}$ solve the following equation with $\tt X_{t_2-t_1}=\tt Y_{t_2-t_1}=0$:
\beq\label{EP3} \beg{cases} \d \tt X_t=  (A \tt X_t+B\tt Y_t-L_1\tt X_t)\d t,\\
\d \tt Y_t=  \big\{Z(\tt X_t,\tt Y_t)-L_2\tt Y_t  \big\}\d t +\si\d W_t,\ \ t\ge t_2-t_1.\end{cases}\end{equation} Then the distribution of $(X_{t_2},Y_{t_2})$ is $\mu_{t_2}$ while that of $(\tt X_{t_2},\tt Y_{t_2})$ is $\mu_{t_1}.$ By the definition of $W$, we have
\beq\label{CA2} W(\mu_{t_1},\mu_{t_2})\le \E (|X_{t_2}-\tt X_{t_2}|+ |Y_{t_2}-\tt Y_{t_2}|).\end{equation}
On the other hand,   \eqref{E1.1}, \eqref{EP3},   $(B2)$ and $(B3)$ imply  that  for any $t\ge t_2-t_1,$
\beg{equation*} \beg{split} &|X_t-\tt X_t|\le \e^{-\ll_1(t-t_2+t_1)} |X_{t_2-t_1}| + \int_{t_2-t_1}^t \e^{-\ll_1(t-s)} (\dd |X_s-\tt X_s| + \|B\|\cdot |Y_s-\tt Y_s|)\d s,\\
 &|Y_t-\tt Y_t|\le \e^{-\ll_1(t-t_2+t_1)} |Y_{t_2-t_1}| + \int_{t_2-t_1}^t \e^{-\ll_1(t-s)} (K_1  |X_s-\tt X_s| + K_2 |Y_s-\tt Y_s|)\d s.\end{split}\end{equation*} Then by \eqref{DEE}, for $t\ge t_2-t_1$
\beq\label{DEE2} \beg{split} &\aa |X_t-\tt X_t|+ |Y_t-\tt Y_t|\\
&\le \e^{-\ll_1(t+t_1-t_2)} (\aa|X_{t_1}|+|Y_{t_1}|) +\ll' \int_{t_2-t_1}^t \e^{-\ll_1(t-s)} (\aa|X_{s}-\tt X_s|+|Y_{s}-\tt Y_s|)\d s.\end{split}\end{equation}   By Gronwall's inequality, we obtain
\beg{align*}  \aa |X_{t_2}-\tt X_{t_2}|+ |Y_{t_2}-\tt Y_{t_2}| &\le (\aa|X_{t_1}|+|Y_{t_1}|)\e^{-\ll_1 t_1}\bigg(1+\ll' \int_{t_2-t_1}^{t_2} \e^{\ll' (t_2-s) }\d s\bigg)\\
  &\le 2(\aa|X_{t_1}|+|Y_{t_1}|)\e^{-(\ll_1-\ll') t_1}.\end{align*}
Since $\sup_{t\ge 0} \E (|X_t|+|Y_t|)<\infty$  due to  Lemma \ref{L4.2}, this together with  \eqref{CA2}  implies \eqref{CA}. The proof is therefore finished.
\end{proof}

\beg{lem} Assume $(B1)$, $(B2)$ and $(B3)$. If    $ \ll_1>\ll'$, then there exists a constant $C>0$ such that for any mild solutions $(X_t,Y_t)$ and $(\tt X_t,\tt Y_t)$ of the equation $\eqref{E1.1}$,
$$|X_t-\tt X_t|+|Y_t-\tt Y_t|\le C(|X_0-\tt X_0|+|Y_0-\tt Y_0|) \e^{-(\ll_1-\ll') t},\ \ t\ge 0.$$\end{lem}

 \beg{proof} Similarly to the proof of \eqref{DEE2}, we have
\beg{align*} &\aa |X_t-\tt X_t|+ |Y_t-\tt Y_t|\\
&\le \e^{-\ll_1t} (\aa|X_{0}-\tt X_0|+|Y_{0}-\tt Y_0|) +\ll' \int_{0}^t \e^{-\ll_1(t-s)} (\aa|X_{s}-\tt X_s|+|Y_{s}-\tt Y_s|)\d s,\ \ t\ge 0.\end{align*} By Gronwall's inequality,
 $$\aa |X_t-\tt X_t|+ |Y_t-\tt Y_t|\le \e^{-(\ll_1-\ll') t } (\aa|X_{0}-\tt X_0|+|Y_{0}-\tt Y_0|),\ \ t\ge 0.$$ This completes the proof. \end{proof}

 \section{Some Examples}

In this section, we present three examples  to illustrate Theorems \ref{T1.1} and \ref{T1.2}, where the first   includes the
 kinetic Fokker-Planck equation discussed in \cite{V} for $V(x)=-\ff 1  2 |x|^2 +\nn W$ with small $\|\nn^2 W\|_\infty$,  the second is highly degenerate in the sense that $m$ can be much larger than $d$,  
  and the last is an infinite-dimensional model.

 \paragraph{Example 5.1.} Let $d=m$ and $\si$ be invertible, $A=0, B=I$, and $Z(x,y)= \nn W(x) -x-y$ for some $W\in C^2(\R^d)$.  If $\|\nn^2 W\|_\infty<1$ is small enough such that 
 \beq\label{EX1} 1>\inf_{r_0\in (0,1)} \Big\{\ff{\|\nn^2W\|_\infty^2}{2r_0(1-\|\nn^2W\|_\infty)(1+\ss{1+4r_0})} +\ff {r_0}2 \big(1+\ss{1+4r_0}\big)\Big\},\end{equation}
   then all assertions in Theorem \ref{T1.1} hold. In particular, \eqref{EX1} holds if $\|\nn^2 W\|_\infty\le \ff 1 2.$   
   
\beg{proof} It is trivial that $(A1)$ and $(A2)$ hold. To verify $(A3)$, let $r>0$ and $r_0\in (0,1)= (0, \|B\|^{-1})$.  By $A=0, B=I$ and the formulation of $Z$, we have 
   \beg{equation*}\beg{split} &\<r^2(x-\bar x) + rr_0 B(y-\bar y), A(x-\bar x)+B(y-\bar y)\>+ \<Z(x,y)-Z(\bar x,\bar y), y-\bar y+rr_0B^*(x-\bar x)\>\\
   &=(r^2-1-rr_0)\<x-\bar x, y-\bar y\> + rr_0|y-\bar y|^2 +\<\nn W(x)-\nn W(\bar x), y-\bar y+rr_0(x-\bar x)\>\\
   &\qquad  - rr_0|x-\bar x|^2 -|y-\bar y|^2.  \end{split}\end{equation*} Take
   \beq\label{DG} r=\ff 1 2 \big(1+\ss{1+4r_0}\big)\end{equation} such that $r^2-1-rr_0=0$, we obtain 
\beg{equation*}\beg{split} &\<r^2(x-\bar x) + rr_0 B(y-\bar y), A(x-\bar x)+B(y-\bar y)\>+ \<Z(x,y)-Z(\bar x,\bar y), y-\bar y+rr_0B^*(x-\bar x)\>\\
&\le -\big(rr_0-\|\nn^2 W\|_\infty rr_0-\gg\big)|x-\bar x|^2- \Big(1- rr_0-\ff{\|\nn W\|_\infty^2}{4\gg}\Big)|y-\bar y|^2,\ \ \gg>0. \end{split}\end{equation*} 
Therefore, $(A_3)$ holds for some constants $r_0\in (0,1)$ and $\theta>0$ if
$$ 1 >\inf_{r_0\in (0,1)} \inf_{\gg\in (0, rr_0-\|\nn^2W\|_\infty rr_0)} \Big(rr_0+\ff{\|\nn W\|_\infty^2}{4\gg}\Big),$$ which is equivalent to \eqref{EX1} due to \eqref{DG}. 
It remains to prove \eqref{EX1} for $\|\nn^2 W\|_\infty\le \ff 1 2.$ Since \eqref{EX1} is trivial for $\|\nn^2W\|_\infty=0$, we assume that $\|\nn^2W\|_\infty\in (0,\ff 1 2].$ In this case we simply take $r_0= \|\nn^2W\|_\infty$ such that
\beg{align*} &\ff{\|\nn^2W\|_\infty^2}{2r_0(1-\|\nn^2W\|_\infty)(1+\ss{1+4r_0})} +\ff {r_0}2 \big(1+\ss{1+4r_0}\big)\\
&< \|\nn^2 W\|_\infty \Big(\ff 1 2+\ff 1 2 \big(1+\ss{3}\big)\Big)\le \ff 1 2 \Big(1+\ff 1 2 \ss 3\Big)<1.\end{align*}
  \end{proof}

 \paragraph{Example 5.2.} Let $\si$ be invertible, $m=kd$ for some natural number $k\ge 2,$ and
 \beg{equation*} \beg{split} & By= (0,\cdots, 0, y)\in \R^{k d},\ \ \ \ y\in \R^d,\\
 &Z(x,y)= b(y)-x_k,\ \ \ \ y\in \R^d, x= (x_1, x_2,\cdots, x_k)\in \R^{kd},\\
 &A(x_1, x_2,\cdots, x_k)= (\gg x_2-x_1, \gg x_3-x_2,  \cdots, \gg x_k- x_{k-1}, 0),\ \ \ \ x_1,\cdots, x_k\in\R^d,\end{split}\end{equation*}
 where $\gg\ne 0$ is a constant, and $b: \R^d\to \R^d$ satisfies 
 \beq\label{C1} |b(y)-b(\bar y)|\le K|y-\bar y|,\ \ \<b(y)-b(\bar y), y-\bar y\>\le -\bb |y-\bar y|^2,\ \ y,\bar y\in \R^d \end{equation} for some constants $K,\bb>0$.   If
 \beq\label{C5}  0<|\gg| < 1\land \ff{2\bb}{2+K^2},\end{equation}  then assertions in Theorem \ref{T1.1} hold.

\beg{proof} It  is easy to see that when $\gg\ne 0$, the rank condition in $(A1)$ holds. Since $b$ is Lipchitz continuous and $\si$ is invertible,
by Theorem \ref{T1.1} it suffices to verify $(A3)$. We simply take $r=1$. For any $r_0\in (0,1)= (0, \|B\|^{-1})$,  we have
 \beg{equation*}\beg{split} &\<r^2(x-\bar x) + rr_0 B(y-\bar y), A(x-\bar x)+B(y-\bar y)\>+ \<Z(x,y)-Z(\bar x,\bar y), y-\bar y+rr_0B^*(x-\bar x)\>\\
 &=  r_0|y-\bar y|^2
 +\sum_{i=1}^{k-1} \Big\{\gg \<x_i-\bar x_i, x_{i+1}-\bar x_{i+1}\>-|x_i-\bar x_i|^2\Big\}
  \\
  &\qquad +\big\<b(y)-b(\bar y), y-\bar y+ r_0(x_k-\bar x_k)\big\>   - r_0  |x_k-\bar x_k|^2 \\
 &\le -(\bb-r_0)|y-\bar y|^2 -r_0|x_k-\bar x_k|^2 + r_0K |y-\bar y|\cdot |x_k-\bar x_k|\\
  &\qquad -\sum_{i=1}^{k-1} \Big\{|x_i-\bar x_i|^2 -\ff{|\gg|} 2  |x_i-\bar x_i|^2- \ff{|\gg|} 2 |x_{i+1}-\bar x_{i+1}|^2\Big\}\\
 &\le -\sum_{i=1}^{k-1} (1-|\gg|)|x_i-\bar x_i|^2-\Big(r-\ff{|\gg|} 2 - \ff{r_0K^2}{4\aa}\Big) |x_k-\bar x_k|^2-(\bb-r_0-\aa r_0) |y-\bar y|^2,\ \ \aa>0.
 \end{split}\end{equation*} So, $(A3)$ holds for some $\theta>0$ provided $|\gg|<1$ and
 $$\sup_{r_0\in (0, 1\land \ff{\bb}{1+\aa}),\aa>0} \Big(r_0-\ff {|\gg|} 2-\ff{K^2r_0}{4\aa}\Big)>0.$$ Letting $r_0\uparrow 1\land \ff\bb{1+\aa}$,
 we conclude that $(A3)$ holds provided $|\gg|<1$ and
 $$\sup_{\aa>0} \Big(1\land\ff{\bb}{1+\aa}\Big)\Big(1- \ff{K^2}{4\aa}\Big)>\ff {|\gg|}2.$$
 By taking $\aa= \ff 1 2 K^2$ we see that this inequality follows from \eqref{C5}.
 \end{proof}

Finally, we  present an example for Theorem \ref{T1.2} in the spirit of
Example 5.2 that $\H_2$ is a subspace of $\H_1.$

\paragraph{Example 5.3.} Let $\{u_i\}_{i\ge 1}$ be an orthonormal basis on   $\H_1$, and let  $\H_2=\overline{\rm span}\{u_{2i}:i\ge 1\}$. Take
    $B=I_{\H_2}$ and    $$ L_1u_{2i}=\ll_iu_{2i},\  L_1u_{2i-1}=\ll_i u_{2i-1},\ \ i\ge 1,$$ where    $0<\ll_i\uparrow\infty$  with
     $\sum_{i\ge 1}\ll_i^{-1}<\infty$. Moreover, let   $L_2=L_1|_{\H_2}$  and
  $$Ax= \gg \ll_1\sum_{i=1}^\infty \<x,u_{2i}\>u_{2i-1},\ \ \ x\in\H_1 $$  for some constant   $\gg\in\R$. Finally, let
  $Z$  satisfy
  $$|Z(x,y)-Z(\bar x, \bar y)|\le \aa \ll_1|x-\bar x|+\bb \ll_1 |y-\bar y| $$  for some constants  $\aa,\bb\ge 0$.
   Then    all assertions in Theorem \ref{T1.1} hold provided
 \beq\label{EEX}  \ss{1+\gg^2} + 4\bb +\ss{(2\bb-1-\ss{1+\gg^2})^2+8\aa}<7.\end{equation}

 \beg{proof} It is easy to see that  $BL_2=L_1 B, AL_1=L_1 A.$ According to  Theorem \ref{T1.2}, it suffices to prove
 \beg{enumerate} \item[(a)] For  some    $\dd>0$   such that  $L_1-A\ge \ll_1-\dd$  and the condition \eqref{C11} hold.
  \item[(b)]  For any   $t_0>0$,
 $Q_{t_0}$ is invertible on  $\H_1$.\end{enumerate}
\paragraph{Proof of  (a)}      We have
  \beg{equation*}\beg{split} &\<(L_1-A)x,x\>   =\<L_2\pi x, \pi x\>-\<Ax,x\>\\
&\ge \ll_1\sum_{i\ge 1} \<x,u_{2i}\>^2 -\gg \sum_{i\ge 1} \<x,u_{2i}\>\<x,u_{2i-1}\>\\
 &\ge (\ll_1-\dd)\sum_{i\ge 1} \<x,u_{2i}\>^2 -\ff{\gg^2}{4\dd}\sum_{i\ge 1}\<x,u_{2i-1}\>^2,\ \ x\in\H_1.\end{split}\end{equation*}
Taking
   $$\dd= \ff{1+\ss{1+\gg^2}}2 \ll_1 $$  such that  $\ff{\gg^2}{4\dd}= \dd-\ll_1$,   we have   $L_1-A\ge \ll_1-\dd$   as required, and
the condition   \eqref{EEX}  is equivalent to  \eqref{C11}.

\paragraph{Proof of   (b)}   We may simply assume $\gg\ll_1=1$, so that
 $$A^* x= \sum_{i=1}^\infty \<x,u_{2i-1}\>u_{2i},\ \ x\in\H_1.$$ Since $A^2=(A^*)^2=0$ and $BB^*$ is the orthogonal projection onto $\H_2$,
 for any $x\in\H_1 $ we have
 \beg{equation*}\beg{split} &\e^{sA}BB^*\e^{sA^*} x = (I+sA) BB^* \{x+ sA^*x\}\\
  &= \sum_{i=1}^\infty\big(\<x,u_{2i}\> +s\<x,u_{2i-1}\>\big)\big\{u_{2i} +su_{2i-1}\}.\end{split}\end{equation*} Then
  \beg{equation*}\beg{split} \<Q_{t_0}x,x\>&= \sum_{i=1}^\infty  \int_0^{t_0} \big\{\<x,u_{2i}\>^2 + 2s \<x,u_{2i-1}\>\<x,u_{2i}\>
  +s^2\<x,u_{2i-1}\>^2\big\}\d s\\
  &= t_0\sum_{i=1}^\infty \Big\{\<x,u_{2i}\>^2 +t_0 \<x,u_{2i-1}\>\<x,u_{2i}\> +\ff{t_0^2}3 \<x, u_{2i-1}\>^2\Big\}\\
  &\ge t_0 \sum_{i=1}^\infty\Big\{(1-r)\<x,u_{2i}\>^2 +\Big(\ff 1 3- \ff 1 {4r}\Big)t_0^2 \<x,u_{2i-1}\>^2\Big\},\ \ r>0.\end{split}\end{equation*}
  Taking $r\in (0,1)$ but
  close enough to $1$, we conclude that $\<Q_{t_0}x,x\>\ge c|x|^2$ holds for some constant $c>0$ and all $x\in \H_1.$ Therefore, $Q_{t_0}$ is invertible.
\end{proof}

\paragraph{Acknowledgement.} The author  would like to thank the referee  for helpful comments.

\beg{thebibliography}{99}

 \bibitem{ATW06} M. Arnaudon, A. Thalmaier, F.-Y. Wang,
  \emph{Harnack inequality and heat kernel estimates
  on manifolds with curvature unbounded below,} Bull. Sci. Math.   130(2006), 223--233.

% \bibitem{BE}  D. Bakry, D. and M. Emery, \emph{Hypercontractivit\'e de
%semi-groupes de diffusion}, C. R. Acad. Sci. Paris. S\'er. I Math.
%299(1984), 775--778.

%\bibitem{ATW14} M. Arnaudon, A. Thalmaier, F.-Y. Wang, \emph{Equivalent Harnack and gradient inequalities for pointwise curvature lower bound,}   Bull. Sci. Math. 138(2014),643--655.

\bibitem{BWY} J. Bao, F.-Y. Wang, C. Yuan, \emph{Hypercontractivity for functional stochastic differential
equations,} Stoch. Proc. Appl. 125(2015) 3636--3656.

 \bibitem{B} F. Baudoin, \emph{ Bakry-Emery meet Villani,} arXiv:1308.4938

 \bibitem{Chen} M.-F. Chen, \emph{From Markov Chains to Non-Equilibrium Particle Systems,} World Scientific, 1992, Singapore.

 \bibitem{DP} G.D.~Prato, J.~Zabczyk, \emph{
Stochastic Equations in Infinite Dimensions,}
Cambridge University Press, 1992.

 \bibitem{DMS} J. Dolbeault, C. Mouhot,  C. Schmeiser, \emph{ Hypocoercivity for kinetic equations with linear relaxation terms,}   C. R. Math. Acad. Sci. Paris 347(2009),  511--516.

 \bibitem{D}  R. Duan, \emph{Hypocoercivity of linear degenerately
dissipative kinetic equations,}  Nonlinearity 24(2011), 2165--2189.

\bibitem{GM} S. Gadat,  L. Miclo,  \emph{Spectral decompositions and $L^2$-operator norms of toy hypocoercive semi-groups,} Kinetic and related models 6(2013), 317--372.

\bibitem{G} L. Gross,  \emph{Logarithmic Sobolev inequalities and contractivity properties of semigroups,}
Lecture Notes in Math. 1563, Springer-Verlag, 1993.

\bibitem{GS} M. Grothaus, P. Stilgenbauer,  \emph{Hypocoercivity for kolmogorov backward evolution equations and applications,} J. Funct. Anal. 267(2014), 3515--3556.

 \bibitem{GW} A. Guillin, F.-Y. Wang,
  \emph{Degenerate Fokker-Planck equations : Bismut formula, gradient estimate  and Harnack inequality,}  J. Diff. Equat. 253(2012), 20--40.

 \bibitem{N}  E. Nelson, \emph{The free Markov field,} J. Funct. Anal. 12 (1973), 211--227.

 \bibitem{S} T. Seidman, \emph{How violent are fast controls?,} Math. Control Signals Systems 1(1988), 89--95.

\bibitem{V} C. Villani, \emph{Hypocoercivity,} Mem. Amer. Math. Soc.  202(950)(2009).

\bibitem{W97} F.-Y. Wang, \emph{Logarithmic Sobolev inequalities on noncompact Riemannian manifolds,} Probab. Theory Relat. Fields
109(1997),  417--424.

\bibitem{W04} F.-Y. Wang, \emph{Spectral gap for hyperbounded operators,} Proc. Amer. Math. Soc. 132(2004), 2629--2638.

\bibitem{W10} F.-Y. Wang, \emph{Harnack inequalities on manifolds with boundary and applications,}    J.
Math. Pures Appl.   94(2010), 304--321.

%\bibitem{W14}   F.-Y. Wang, \emph{Derivative formula and gradient estimates for  Gruschin type semigroups,}  J. Theo. Probab. 27(2014), 80--95.

\bibitem{Wbook} F.-Y. Wang, \emph{Harnack Inequalities and Applications for Stochastic Partial Differential Equations,} Springer, 2013, Berlin.

%\bibitem{W14} F.-Y. Wang, \emph{Criteria on spectral gap of Markov operators,} J. Funct. Anal. 266(2014), 2137--2152.

\bibitem{WY11} F.-Y. Wang, C. Yuan, \emph{Harnack inequalities for functional SDEs with multiplicative noise and applications,}   Stoch. Proc. Appl.    121(2011), 2692--2710.

 \bibitem{WZ} F.-Y. Wang, T. Zhang, \emph{Gradient estimates for stochastic evolution equations with non-Lipschitz coefficients,}  J. Math. Anal. Appl. 365(2010), 1--11.

 \bibitem{WZ2} F.-Y. Wang, T. Zhang, \emph{Degenerate SDEs in Hilbert spaces with rough drifts, }  Infin. Dimens. Anal. Quant. Probab. Relat. Top.   18 (2015), no. 4, 1550026, 25 pp

\bibitem{Wu} L. Wu, \emph{ Uniformly integrable operators and large deviations for Markov processes,}
J. Funct. Anal. 172 (2000), 301--376.

\end{thebibliography}
\end{document}